\documentclass[12pt]{article}
\usepackage{amsmath}
\usepackage{graphicx,psfrag,epsf}
\usepackage{enumerate}
\usepackage{natbib}
\usepackage{url} 
\usepackage{amssymb}
\usepackage{changebar}

\newcommand{\blind}{1}

\newcommand{\comment}[1]{}

\newfont{\cmrten}{cmr10}
\newtheorem{theo}{Theorem}[section]
\newtheorem{lem}{Lemma}[section]

\newtheorem{ex}{Example}
\newtheorem{definition}{Definition}[section]

\newtheorem{remark}{Remark}[section]

\newcommand{\vecc}[1]{\mbox{\boldmath $#1$}}

\newcommand{\Real}{\mbox{\rm I\kern-.23em\hbox{R}}}
\newcommand{\commentt}[1]{}

\def\ta{\tilde{a}}
\def\tb{\tilde{b}}
\def\P{m_{\theta^{*}}}
\def\X{X}
\def\x{x}

\def\by{\vecc{t}}

\def\bz{\vecc{z}}
\def\bZ{\vecc{Z}}
\def\b1{\vecc{1}}
\def\cp{\P^{(c)}}

\def\mhat{ m_{\hat\theta}}
\def\cmhat{ m^{(c)}_{\hat\theta}}
\def\cmhatw{ m^{(wc)}_{\hat\theta}}

\def\rsqp{\rho^2_{\P}}
\def\lsqp{\lambda^2_{\P}}

\def\aw{\hat{a}^{(w)}}
\def\bw{\hat{b}^{(w)}}
\def\yw{\bar{T}^{(w)}}
\def\mw{\bar{m}^{(w)}_{\hat\theta}}

\def\bold1{{\bf 1}}

\def\veccz0{\vecc{z}_0}

\def\convinprob{\buildrel P \over \longrightarrow}

\def\convindist{\buildrel d \over \longrightarrow}

\def\convinprob{\buildrel P \over \longrightarrow}

\def\convindist{\buildrel d \over \longrightarrow}

\begin{document}


\def\spacingset#1{\renewcommand{\baselinestretch}%
{#1}\small\normalsize} \spacingset{1}

\include{titlepage_blind}


\if1\blind
{

\title{Prediction Accuracy  Measures for a Nonlinear Model and for Right-Censored Time-to-Event  Data}

\author{
Gang Li\footnote{Gang Li is Professor of  Biostatistics and Biomathematics,
University of California, Los Angeles, CA
90095-1772, USA. (E-mail: vli@ucla.edu).
The research of Gang Li was partly supported by National Institute of Health Grants P30 CA-16042, UL1TR000124-02, and P01AT003960}  and Xiaoyan Wang\footnote{Xiaoyan Wang is Adjunct Assistant Professor, Division of General Internal Medicine and Health Services Research,
University of California, Los Angeles, CA
90095-1736, USA. (E-mail:  xywang@mednet.ucla.edu).
}
\\[1.0cm]
{\sl Running Title:  Prediction Accuracy Measures}\\
}

\date{}

\maketitle

} \fi

\if0\blind
{
  \bigskip
  \bigskip
  \bigskip
  \begin{center}
    {\LARGE\bf Prediction Accuracy Measures for a Nonlinear Model and for an Event Time Model with Right-Censored Data}
\end{center}
  \medskip
} \fi

\bigskip
\begin{abstract}
This paper studies prediction summary measures for a prediction function under a general setting in which the model is allowed to be misspecified and the prediction function is not required to be the conditional mean response. We show that the $R^2$ measure based on a variance decomposition is insufficient to summarize the predictive power of a nonlinear prediction function. By deriving a prediction error decomposition, we introduce an additional measure,  $L^2$, to augment the $R^2$ measure. When used together, the two measures provide a complete summary of the predictive power of a prediction function.   Furthermore, we extend these measures to right-censored time-to-event data by establishing right-censored data analogs of the variance and prediction error decompositions. We illustrate the usefulness of the proposed measures with simulations and real data examples.  Supplementary materials for this article are available online.
\end{abstract}

\noindent%
\vfill

\newpage
\spacingset{1.45} 


\section{Introduction}\label{section1}

In this paper we develop prediction accuracy ymeasures for  a nonlinear model and for right-censored time-to-event data. In addition to evaluating a model's prediction performance, prediction Accuracymeasures are useful for assessing the practical importance of predictors and for comparing competing models that are not necessarily nested nor correctly specified.

By far, the most commonly used prediction accuracy measure for a linear model is the
R-squared statistic, or coefficient of determination.
Let $Y$ be a real-valued random variable and $\X$ be a vector of $p$ real-valued explanatory random variables or covariates. Assume that one observes a random sample $(Y_1, \X_1), \ldots, (Y_n, \X_n)$
from the distribution of $(Y,\X)$. The R-squared statistic is defined as
\begin{equation}
\label{rsquarelinear}
R^2
 =1 - \frac{\sum_{i=1}^n (Y_i-\hat Y_i )^2}{\sum_{i=1}^n (Y_i-\bar Y)^2},
\end{equation}
where  $\hat Y_i=a + b^T X_i$ is the least squares predicted value for subject $i$.
The $R^2$ statistic has the straightforward interpretation as the proportion of variation of $Y$ which is explained by the least squares prediction function due to the  following decomposition:
\begin{eqnarray}
\label{decomposition}
\sum_{i=1}^n (Y_i-\bar Y)^2&=& \sum_{i=1}^n (\hat Y_i -\bar Y)^2 + \sum_{i=1}^n (Y_i-\hat Y_i)^2. \\
\mbox{total variation}&=& \mbox{explained variation } +\mbox{ unexplained variation } \nonumber
\end{eqnarray}
Despite its popularity in linear regression,  the  $R^2$ statistic defined by (\ref{rsquarelinear})  is not readily applicable to  a nonlinear model since the decomposition (\ref{decomposition}) no longer holds.  In the past decades, much efforts have been devoted to extending the R-squared statistic to nonlinear models. Among others,  the pseuodo $R^2$ statistics
for a nonlinear model include likelihood-based measures \citep{goodman1971analysis, mcfadden1973conditional, maddala1986limited, cox1989analysis, magee1990r, nagelkerke1991note}, information-based measures \citep{mcfadden1973conditional, kent1983information}, 
ranking-based measures \citep{harrell1982evaluating}, variation-based measures \citep{theil1970estimation, efron1978regression,haberman1982analysis, hilden1991area, cox1992comment, ash1999r2},
and the multiple correlation coefficient measure \citep{mittlbock1996explained, zheng2000summarizing}.
However, none of the existing pseudo $R^2$ measures are motivated directly from a variance decomposition and none have received the same widespread acceptance
as the classical $R^2$ for linear regression.
Interested readers are  referred to \citet{zheng2000summarizing} for an excellent survey of existing pseudo $R^2$ measures and further references on this topic.  

The first goal of this paper is to develop prediction accuracy measures for a prediction function under a general setting in which the model is allowed to be misspecified and the prediction function may be different from the conditional expected response.  We begin with defining population prediction accuracy measures. Based on a simple variance decomposition, we define a $\rho^2$ measure as the proportion of the explained variance of $Y$ by a corrected prediction function.  It can be shown that the $\rho^2$ parameter is identical to the squared multiple correlation coefficient between the response and the predicted response.
Since it describes the proportion of the explained variance by the corrected prediction function, which in general is not the same as the uncorrected prediction functions, the squared multiple correlation coefficient, a popular pseudo $R^2$,  is not sufficient to summarize the predictive power of nonlinear models.  As a remedy, we derive another parameter, named $\lambda^2$, as the proportion of the explained prediction error by the corrected prediction function based on
a mean-squared prediction error decomposition.  The parameter $\lambda^2$ measures how close the
uncorrected prediction function is to its corrected version. The two parameters characterize complementary aspects regarding the predictive accuracy of the prediction function. When used in combination, they provide a complete summary of the predictive power of the uncorrected prediction function. We further obtain finite sample versions of the variance and prediction error decompositions, define the corresponding sample prediction accuracy measures, namely $R^2$ and $L^2$, and establish their asymptotic properties.  It is worth noting that for the least squares prediction function under the linear model, the $L^2$ measure degenerates to 1 and therefore only $R^2$ is needed to describe its predictive power in the classical linear regression analysis.

The second goal of the paper is to develop new prediction accuracy measures for an event time model based on right censored time-to-event data. Note that it is challenging to extend the $R^2$ definition (\ref{rsquarelinear})  to right-censored data even for the linear model. A variety of  pseudo $R^2$ measures and other loss functions have been proposed for event time models with right-censored data \citep{kent1988measures, korn1990measures, graf1999assessment, schemper2000predictive,royston2004new, o2005explained, stare2011measure}. For example, 
the EV option in the SAS PHREG procedure gives a generalized $R^2$ measure proposed by 
\citet{schemper2000predictive} for  Cox's (1972)  proportional hazards model. A more recent proposal by \citet{stare2011measure} uses explained rank information, which is applicable to a wide range of event time models.  \citet{stare2011measure} also gave a thorough literature review of prediction accuracy measures  for event time models. 
 We highlight that for linear regression,
none of the existing pseudo $R^2$  measures for right censored data reduce to the classical R-squared statistic  in the absence of censoring. Moreover, under a correctly specified model, they do not converge to the nonparametric population R-squared value
$\rho_{NP}^2\equiv var{(E(Y|X))}/var(Y)$, the proportion of the explained variance by $E(Y|X)$, as the sample size grows large.  Finally, as shown in Section 4 (Table 1) that the pseudo $R^2$ measures of 
\citet{schemper2000predictive, stare2011measure} are not suitable for comparing unnested Cox's models with possibly different baseline hazards and could remain constant when the nonparametric population R-squared value
$\rho_{NP}^2$ varies from 0 to 1. In this paper, we derive a 
variance and a prediction error decomposition for right censored data. These decompositions allows us to define a pair of  prediction accuracy measures, $R^2$ and $L^2$,  for an event time model with right-censored data in exactly the same way as uncensored data.  The proposed measures possess many  appealing properties that most existing pseudo $R^2$ measures do not have. First, for the linear model with no censoring,  our $R^2$ statistic reduces to the classical coefficient of determination and $L^2$ reduces to 1. Second,  when  the prediction is  the conditional mean response  based on a correctly specified model,  our $R^2$ statistic is a consistent estimate of the nonparametric coefficient of determination $\rho_{NP}^2$, and $L^2$ converges to 1 as the sample size grows large. Third, our method is applicable to any event time model with right-censored data. Fourth,  our measures are defined without requiring the model to be correctly  specified. Lastly,  our measures can be used to compare unnested models.

The rest of the paper is organized as follows. In Section 2.1, we define a pair of population prediction accuracy measures for  a general prediction function from a possibly mis-specified model by deriving a variance decomposition and a  mean squared prediction error decomposition.  Sample measures based on independent and identically distributed complete data are then proposed and studied in Section 2.2. Section 3 discusses how to extend these measures to event time models with right-censored data. Section 4 presents simulation studies to illustrate the performance of the proposed sample measures and compare them with some existing measures in the literature. Real data illustrations are given in Section 5. Proofs of theoretical results are deferred to Appendix. Final remarks are provided in Section 6.

\section{Prediction Summary Measures for a Nonlinear Model}

Denote by $F(y|x)=P(Y\le y|X=x)$ and $\mu(x) = E(Y|\X=\x)$ the true conditional distribution function and the true conditional expectation of $Y$ given $X=x$, respectively.

Consider a regression model of $Y$ on $\X$ described by a family of conditional distribution functions
${\cal M }=\{ F_{\theta}(y|\x): \theta\in \Theta\}$, where the parameter $\theta$ is either finite dimensional or infinite dimensional.  For example,
$F_{\theta}(y| \x) = \Phi((y-\alpha -\beta^T \x)/\sigma)$ for the linear regression model with a normal $N(0,\sigma^2)$ random error, where
$\theta = (\alpha,\beta^T, \sigma^2)$ and $\Phi$ is the standard normal cumulative distribution function. The
Cox (1972) proportional hazards model is an example of a semi-parametric regression model
with
$F_{\theta}(y| \x) = 1 - \{ 1-F_0(y)\}^{\exp(\beta^T\x)}$
where $\theta =(\beta, F_0)$ consists of a finite dimensional regression parameter
$\beta$ and an infinite dimensional unknown baseline distribution function $F_0$. We allow the model
${\cal M}$ to be misspecified in the sense that ${\cal M}$ may not include the true conditional distribution function $F(y|x)$ as a member.

For any $\theta\in \Theta$, let $m_{\theta}(\X)$ be a prediction function of $Y$ obtained as a
functional of  $F_{\theta}(\cdot|\X)$.  Common examples of  $m_{\theta}(\X)$  include the conditional mean response defined by  $m_{\theta}(\x) = \int y dF_{\theta}(y|\x)$ and
the conditional median response $m_{\theta}(\x) =F^{-1}_{\theta}(0.5|\x)$.
Assume that $\hat\theta$ is a sample statistic such that
as $n\to\infty$,
\begin{equation}
\label{C1}
\hat\theta \convinprob \theta^*,   \quad \mbox{for some $\theta^* \in  \Theta$. }
\end{equation}
For example,  if $\hat\theta$ is  the maximum likelihood estimate for a parametric model, then
under some regularity conditions
$\hat\theta$  converges in probability to a well-defined limit, $\theta^*$, even when the model is misspecified
\citep{huber1967behavior}. If the
model is correctly specified, then $\theta^*$ is the true parameter value. On the other hand, if the model is misspecified, then $\theta^*$ is the parameter that minimizes the Kullback-Leibler Information Criterion
\citep{akaike1998information}.

In this section, we first develop population prediction accuracy measures for $m_{\theta^*}(\X)$, which can be regarded as the asymptotic accuracy measures for the predictive power of $m_{\hat\theta}(\x)$.  Sample prediction accuracy measures for  $m_{\hat\theta}(\X)$ are then derived
accordingly and their asymptotic properties are studied.

\subsection{Population Prediction Summary Measures}

For any  $p$-variate function $P(x)$, define
$$
MSPE(P(X)) = E\{Y-P(\X)\}^2.
$$
as the  {\it mean squared prediction error } ($MSPE$) of  $P(\X)$ for predicting $Y$.

In general, one would expect a good prediction function $P(X)$ of $Y$ to possess at least the following basic
properties:  i)  $E\{P(X)\} =\mu_Y$, and ii) $MSPE(P(X))\le MSPE(\mu_Y)$, where $\mu_Y =E(Y)$ is the best prediction among all constant (non-informative) predictions of $Y$ as measured by $MSPE$. However,
such minimal requirements  are not always satisfied by  $\P(\X)$ when the model
${\cal M} $ is possibly misspecified or when the prediction is not based on the conditional mean response.  Below we introduce a linear correction of $\P(\X)$ so that the corrected prediction function always satisfies
these minimal requirements.
\begin{definition}
\label{definition1}
The linearly corrected prediction function of $\P(\X)$ is defined as
\begin{equation}
\label{cp}
\cp(\X) = \mu_Y + \frac{cov(Y, \P(\X))}{var(\P(\X))} [\P(\X) - E\{ \P(X)\}].
\end{equation}
\end{definition}
It is straightforward to show that $\cp(X)$ has the following  properties.
\begin{enumerate}

\item[(i)]
$\cp(X) = \ta + \tb\P(X)$, where
$
(\tilde{a}, \tilde{b})=\arg \min_{\alpha, \beta} E\{ Y -  (\alpha +\beta \P(\X))\}^2 $;
\item[(ii)] $E(\cp(\X)) = \mu_Y$;

\item[(iii)]
 $MPSE(\cp(\X)) \le MPSE(\mu_Y)$;
\item[(iv)]
$MPSE(\cp(\X)) \le MPSE(\P(X))$.
\end{enumerate}
It follows from (i) and (ii) that $\cp(\X)$ is the best unbiased prediction of $Y$ among all linear functions of $\P(\X)$.
Moreover, the corrected function facilitates two elementary decompositions as stated in Lemma ~\ref{prop1} below.

\begin{lem}
\label{prop1}
Let $\cp(X)$ be the corrected prediction function of $\P(X)$ defined by  (\ref{cp}). Then,
\begin{enumerate}
\item[(a)] (Variance decomposition)
\begin{eqnarray}
\label{decomposition1}
var(Y)&=& E\{\cp(\X) -\mu_Y\}^2 +E\{Y - \cp(\X)\}^2,\\
 &=& \mbox{explained variance} +  \mbox{unexplained variance} \nonumber
\end{eqnarray}
where the first and second terms on the right hand side represent respectively the explained variance and the unexplained  variance of $Y$ by $\cp(\X)$.
\item[(b)] (Prediction Error Decomposition)
\begin{eqnarray}
\label{decomposition2}
MSPE(\P(X))&=& E\{Y-\cp(X)\}^2 + E\{\cp(X) - \P(X)\}^2\\
&=& \mbox{explained prediction error}
+\mbox{ unexplained prediction
error } \nonumber
\end{eqnarray}
where the first and second terms on the right hand side can be interpreted as the explained prediction error and unexplained prediction error of $\P(\X)$ by $\cp(\X)$.
\end{enumerate}

\end{lem}

Based on the above decompositions, we introduce the following prediction accuracy measures.
\begin{definition}
\label{definition2}
Define
\begin{eqnarray}
\label{rsquare}
\rsqp&=&1- \frac{E\{Y - \cp(X)\}^2}{var(Y)} = \frac{E\{\cp(X) -\mu_Y\}^2}{var(Y)},
\end{eqnarray}
to be the proportion of the variance of $Y$ that is explained by $\cp(\X)$, and
\begin{eqnarray}
\label{lsquare}
\lambda^2_{\P}&=& \frac{MSPE(\cp(X)) \}^2}{MSPE(\P(X))}=1 -  \frac{E\{\cp(X) -\P(\X)\}^2}{MSPE(\P(X))} .
\end{eqnarray}
to be the proportion of the $MSPE$ of $\P(X)$ that is explained by
 $\cp(X)$.
\end{definition}

\begin{remark}
\label{interpretation}
{\rm The parameters $\rsqp$ and $\lambda^2_{\P}$ measure two distinct, yet complementary aspects regarding the prediction accuracy of $\P(\X)$:  $\rsqp$  measures the predictive power of the corrected prediction function $\cp(\X)$, whereas
 $\lambda^2_{\P}$  measures how close $\P(\X)$ is to $\cp(\X)$. When used together, they provide a complete accuracy of the predictive power of 
 the uncorrected prediction function $\P(X)$. Note that
$0\le \rsqp \le 1$ and $0\le  \lambda_{\P}^2 \le 1$.  Moreover,  $\rsqp=1$ and $ \lsqp=1$ if and only if $\P(X)=Y$ with probability 1.   So $\P(X)$ has high predictive power if both measures are close to 1.    If
$\rsqp$ is large, but  $\lambda^2_{\P}$ is small, then $\P(X)$ does not have good predictive power even though the corrected prediction $\cp(X)$ does. Lastly,  if $\rsqp$ is small,
 then $\cp(X)$ and consequently $\P(X)$ both do not have good prediction power regardless the magnitude of $\lambda^2_{\P}$.}
\end{remark}

\begin{remark} (Geometric Interpretation).
\label{remark2.1}
{\rm  One may gain more insight about these parameters by examining the geometric relationship between the related quantities.
Define the $L_2$-distance between any two real-valued random variables $\xi$ and
$\eta$ by
$
d_2(\xi,\eta)= \left\{E(\xi -\eta )^2 \right\}^{\frac{1}{2}}.
$
The geometric relationship between $Y$, $\mu_Y$, $\P(\X)$, $\cp(\X)$, and $\mu(\X)$ are  depicted in Figure 1, in which ${\cal P}(X)$ denotes the space of all real-valued functions of $X$.
\begin{center}
[Insert Figure 1 approximately here]
\comment{
\begin{figure}[htbp]
\begin{center}
\includegraphics[scale=0.60]{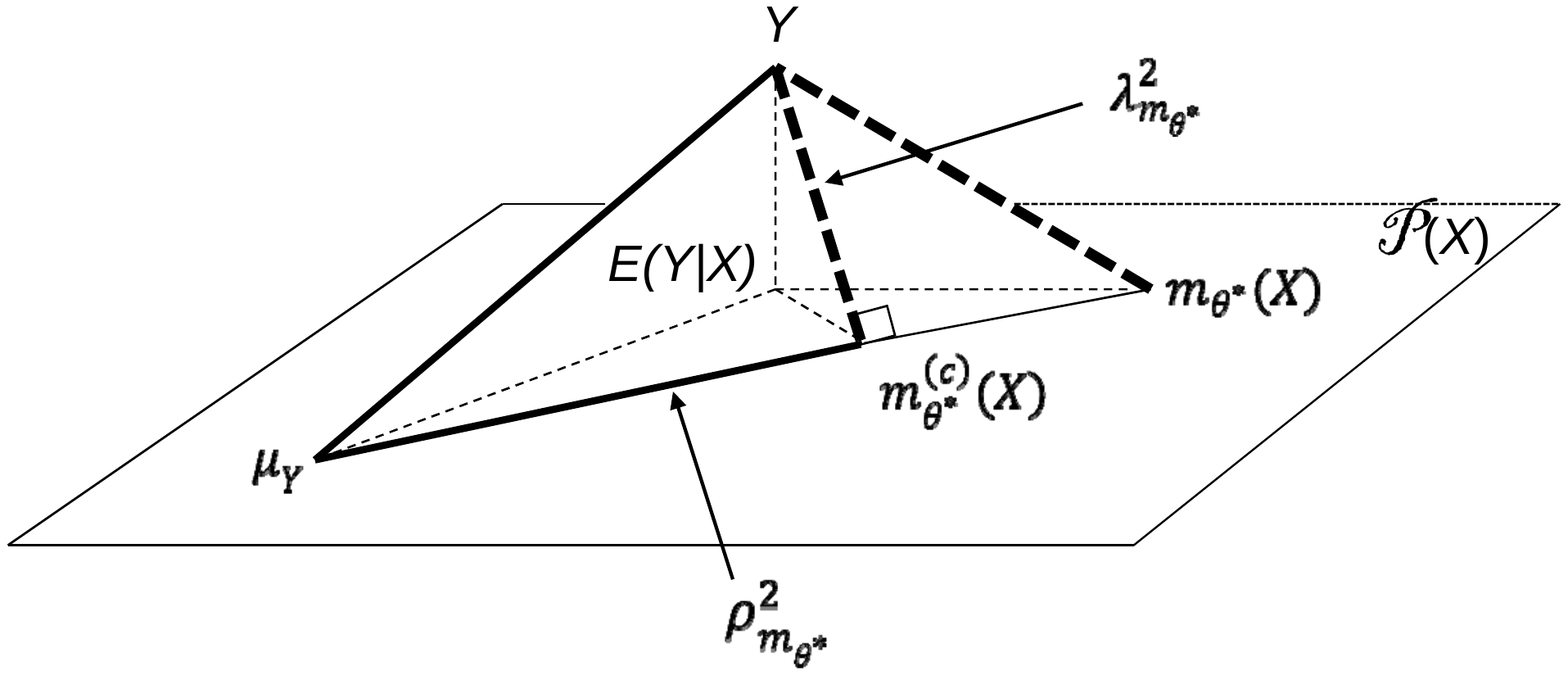}
\caption{\label{fig1} Geometric interpretation of $\rsqp$ and $\lambda^2_{\P}$}
\end{center}
\end{figure}
}
\end{center}
As illustrated in Figure 1,  $\cp(\X)$ is the projection of $Y$ onto the subspace of all linear functions of $\P(X)$ and $\mu(X)$ is the projection of $Y$ onto ${\cal P}(X)$. The variance decomposition in Lemma \ref{prop1}(a) corresponds to the Pythagorean theorem for the triangle
$(Y,\cp(X), \mu_Y)$ that leads to the definition of $\rsqp$. The prediction error decomposition is the Pythagorean theorem for the triangle
$(Y,\cp(X), \P(X))$ that defines $\lambda^2_{\P}$.  }
\end{remark}

\begin{remark}
\label{remark2.2}
(Interpretation of $\lambda_{\P}^2$ as a measure of the prediction bias  for the mean regression function $\mu(X)$).
{\rm
Assume that $\P(X)$ is a nonlinear prediction function. It is easily seen  that if $\P(X)=\mu(X)$, then $\lsqp=1$.  Thus, $\lambda^2_{\P}<1$ implies that $\P(X)\ne \mu(X)$.  In particular, if $\P(X)$ is the conditional mean response under model ${\cal M}$, then $\lambda^2_{\P}<1$ implies that the model  is mis-specified.}
\end{remark}

It is also seen from Figure 1 that the Pythagorean theorem for
the triangle $(Y, \mu(X), \mu_Y)$  corresponds to the well known variance decomposition
\begin{eqnarray}
\label{nonparametric-decomposition}
var(Y) &=& var(\mu(X)) + E(var(Y|X))\nonumber\\
&=& \mbox{explained variance by $\mu(X)$} +\mbox{ unexplained variance}. \nonumber
\end{eqnarray}
We refer the proportion of explained variance by $\mu(X)$:
\begin{equation}
\label{nonparametricR}
\rho^2_{NP} \equiv 1 - \frac{E(Y-\mu(X))^2}{var(Y)} = \frac{var(\mu(X))}{var(Y)},
\end{equation}
 as   the
{\it nonparametric coefficient of determination}. Note that $\rho_{NP}$ is the ``correlation ratio'' studied previously by \citet{renyi1959measures}.

The next theorem summarizes some fundamental  properties of $\rsqp$ and $\lsqp$.
\begin{theo}
\label{prop2}

\begin{enumerate}
\item[(a)]
 Let  $\rho(\xi, \eta)$ denote the correlation coefficient between  two random variables $\xi$ and $\eta$. Then, $\rsqp = [\rho (Y, \P(X))]^2$;

\item[(b)]
(Linear Prediction).
Let $BLUE(X)=a+b^T X$ be
the best linear unbiased estimator (BLUE) of $Y$,
where $(a,b)=\arg\min_{\alpha,\beta} E\{Y-(\alpha+\beta^T X)\}^2$. Then (i)
$BLUE^{(c)}(X)= BLUE(X)$; (ii)
${\lambda^2_{BLUE}}\equiv 1$; (iii) $\rho^2_{BLUE} $ is equal to the population value of the classical coefficient of determination for linear regression.
\item[(c)]
If $\P(X)=\mu(X)$, then
$\lambda_{\P}^2\equiv 1$, and
$\rho^2_{\P}=\rho^2_{NP}$, where $\rho^2_{NP}$ is the
nonparametric coefficient of determination defined by (\ref{nonparametricR});

\item[(d)]
(Maximal $\rho^2$). Let $\rho^2_{NP}$ be defined by (\ref{nonparametricR}). Then
$$\rho_{NP}^2 = \max_{Q\in {\cal P}(X)}\{\rho^2_Q\}$$
where ${\cal P}(X)$ is the space of all $p$-variate functions $Q(X)$ of $X$.
In other words,  $\rho^2_{NP}$ is the {\it maximal coefficient of determination} over all prediction functions $Q(X)$.

\end{enumerate}

\end{theo}

\comment{
\begin{remark}
\label{classification}
(Classification Accuracy Measures).  Although the accuracy measures $\rho^2$ and $\lambda^2$ are derived in the context of the prediction problem, the decompositions in Lemma \ref{prop1} do not require the response $Y$ to be continuous.  An important example of a non-continuous outcome is the binary response $Y$ modeled by a  logistic regression model $\log\frac{m_{\theta}(X)}{1-m_{\theta}(X)} = \theta^{'} X$. In a classification problem, the mean response function $\P(X)$ is used to classify a subject into 0 or 1 with respect to a cutoff value.  The $\rsqp$ measure offers  an appealing alternative to the commonly used ROC method to summarize the  power of $\P(X)$ as a marker for discriminating between 0 and 1 without reference to a specific cutoff point.

\end{remark}
}

\subsection{Sample Prediction Summary Measures}
\label{sect2.2}
Assume that one observes a random sample $(Y_1,\X_1), \ldots, (Y_n, \X_n)$
 of $n$ independent and identically distributed (i.i.d.) replicates of $(Y,\X)$.
 Now we derive sample accuracy measures for the predictive power of $m_{\hat\theta}(X)$,
where  $\hat\theta=\hat\theta(Y_1,\X_1,\ldots, Y_n, \X_n)$ is a sample statistic satisfying (\ref{C1}).

We first give a finite sample version of the decompositions in Lemma \ref{prop1}.
\begin{lem}
\label{prop3}
Define
\begin{equation}
\label{Q}
\cmhat (\x)=\hat a + \hat b \mhat(\x),
\end{equation}
to be the linearly corrected function for $\mhat(x)$,
where
$\hat a = \bar Y - \hat b \bar m_{\hat\theta},
\hat b = \frac{\sum_{i=1}^n (Y_i-\bar Y)\{m_{\hat\theta}(\X_i) -\bar m_{\hat\theta}\}}{\sum_{i=1}^n \{ m_{\hat\theta}(\X_i) -\bar m_{\hat\theta}\}^2},
$
$
{\bar Y} = n^{-1} \sum_{i=1}^n Y_i,
$
 and $\bar m_{\hat\theta}=n^{-1} \sum_{i=1}^n m_{\hat\theta}(\X_i)$.
In other words, $\cmhat(\x)$ is the ordinary least squares regression function obtained by linearly regressing $Y_1, \ldots, Y_n$ on
$m_{\hat\theta}(\X_1), \dots, m_{\hat\theta}(\X_n)$.
Then
\begin{enumerate}
\item[(a)] (Variance Decomposition)
\begin{equation}
\label{sample-variance-decomp}
\sum_{i=1}^n (Y_i -{\bar Y})^2 = \sum_{i=1}^n (\cmhat(X_i)- \bar Y)^2
+ \sum_{i=1}^n (Y_i - \cmhat (X_i))^2;
\end{equation}

\item[(b)] (Prediction Error Decomposition)
\begin{equation}
\label{sample-pred-error-decomp}
\sum_{i=1}^n (Y_i - \mhat (X_i) )^2 = \sum_{i=1}^n (Y_i -\cmhat(X_i))^2
+\sum_{i=1}^n (\cmhat(X_i)- \mhat(X_i))^2.
\end{equation}
\end{enumerate}
\end{lem}

The sample version of $\rho^2$ and $\lambda^2$ are then defined by
\begin{equation}
\label{Rsquare}
R^2_{\mhat} =  \frac{\sum_{i=1}^n (\cmhat(X_i)- \bar Y)^2}{\sum_{i=1}^n (Y_i -{\bar Y})^2},
\end{equation}
and
\begin{equation}
\label{Lsquare}
L^2_{\mhat}=\frac{\sum_{i=1}^n (Y_i -\cmhat(X_i))^2}{\sum_{i=1}^n (Y_i - \mhat (X_i) )^2},
\end{equation}
where $R^2_{\mhat} $ is  the proportion of variation of $Y$ explained by $\cmhat(X)$ and $L^2_{\mhat}$ is the proportion of prediction error of $\mhat(X)$ explained by $\cmhat(X)$.

\begin{remark}
\label{sample-mesures}
Similar to Theorem \ref{prop2}(a),
$
R^2_{\mhat} = \{r(Y,m_{\hat\theta}(\X))\}^2 $
where $r(Y,m_{\hat\theta}(\X))$ is the Pearson correlation coefficient between $Y$ and $m_{\hat\theta}(\X))$.
It can also be easily verified that if
$\mhat(x)$ is the fitted least squares regression line from a linear model, then $L^2_{\mhat} \equiv 1$
and $R^2_{\mhat} $ is identical to the classical coefficient determination for the linear model.
\end{remark}

Below we give the asymptotic properties of  $R^2_{\mhat} $ and $L^2_{\mhat}$.
\begin{theo}
\label{asymptotics}
Assume condition (\ref{C1}) holds.   Assume further that $\P(x)$ is a bounded function and
\begin{equation}
\label{C2}
\mhat(x) \convinprob  \P(x) \quad \mbox{uniformly in $x$}.
\end{equation}
Then, as $n\to\infty$,
\begin{enumerate}
\item[(a)] (Consistency)
$$
R^2_{\mhat}  \convinprob \rho^2_{\P}, \quad\mbox{and}\quad L^2_{\mhat}  \convinprob \lambda^2_{\P};
$$
\item[(b)]
(Asymptotic normality)
\begin{eqnarray*}
\sqrt{n} (R^2_{\mhat}  -\rsqp) \convindist N(0, \sigma_{\rho}^2),
\quad\mbox{and}\quad
\sqrt{n} (L^2_{\mhat}  -\lsqp) \convindist N(0, \sigma_{\lambda}^2),
\end{eqnarray*}
where $\sigma_{\rho}^2$ and $\sigma_{\lambda}^2$ are the asymptotic variances.
\end{enumerate}
\end{theo}

The asymptotic results allow one to assess the variability of the sample measures  $R^2_{\mhat} $ and $L^2_{\mhat}$ and  obtain confidence interval estimates for the corresponding population parameters. In practice, the
bootstrap method \citep{efron1994introduction} or a transformation-based method would be
more appealing than the normal approximation method because the sampling distributions of  $R^2_{\mhat}$ and $L^2_{\mhat}$ can be skewed, especially near 0 and 1.

\section{Sample Prediction Summary Measures for Right Censored Data}

In this section we extend the prediction accuracy measures $R^2_{\mhat}$ and $L^2_{\mhat}$ developed in the previous section to an event time model with right censored time-to-event data.  Recall that we consider a regression
model of $Y$ on $\X$ described by a family of conditional distribution functions
${\cal M }=\{ F_{\theta}(y|\x): \theta\in \Theta\}$, where the parameter $\theta$ is either finite dimensional or infinite dimensional.
Let $T=\min\{Y,C\}$ and $\delta=I(Y\le C)$, where $C$ is an censoring random variable that is assumed to be
independent of $Y$ given $X$.
Assume that one observes a right censored sample of $n$ independent and identically distributed triplets
$(T_1,\delta_1,X_1),\ldots, (T_n, \delta_n, X_n)$ from the distribution of $(T,\delta, X)$.

Assume that $\hat\theta=\hat\theta (T_1,\delta_1,X_1,\ldots, T_n, \delta_n, X_n)$ is a sample statistic satisfying
(\ref{C1}).
Apparently the sample prediction accuracy measures defined in (\ref{Rsquare}) and (\ref{Lsquare}) are no longer applicable
to right censored data because $Y$ is not observed for everything subject.  Below we obtain right-censored data analogs of the uncensored data decompositions
(\ref{sample-variance-decomp}) and (\ref{sample-pred-error-decomp}), and define prediction summary
measures for right censored data.

\begin{lem}
\label{prop5}
Let $w_1,\ldots, w_n$ be a set of nonnegative real numbers satisfying $\sum_{i=1}^n w_i =1$
Define
\begin{equation}
\label{Q1}
\cmhatw (\x)=\aw + \bw \mhat(\x),
\end{equation}
to be a linearly corrected function for $\mhat(x)$,
where
$\aw  = \yw - \bw \mw$,
$\yw =  \sum_{i=1}^n w_i T_i$,
$
\bw = \frac{\sum_{i=1}^n  w_i(T_i-\yw)\{m_{\hat\theta}(\X_i) -\mw\}}{\sum_{i=1}^n w_i \{ m_{\hat\theta}(\X_i) -\mw\}^2}
$, and
$\mw = \sum_{i=1}^n w_im_{\hat\theta}(\X_i)$. In other words, $\cmhatw(\x)$ is the fitted  regression function from the weighted least squares linear regression of  $Y_1, \ldots, Y_n$ on
$m_{\hat\theta}(\X_1), \dots, m_{\hat\theta}(\X_n)$ with weight $W=diag\{w_1,\ldots,w_n\}$.
Then
\begin{enumerate}
\item[(a)] (Weighted Variance Decomposition for $T$)
\begin{equation}
\label{sample-variance-decomp2}
\sum_{i=1}^n w_i\{T_i -\yw\}^2 = \sum_{i=1}^n w_i \{\cmhatw(X_i)- \yw\}^2
+ \sum_{i=1}^n w_i \{T_i - \cmhatw (X_i)\}^2;
\end{equation}

\item[(b)] (Weighted Prediction Error Decomposition for $T$)
\begin{equation}
\label{sample-pred-error-decomp2}
\sum_{i=1}^n w_i\{T_i - \mhat (X_i) \}^2 = \sum_{i=1}^n w_i\{T_i -\cmhatw(X_i)\}^2
+\sum_{i=1}^n w_i \{\cmhatw(X_i)- \mhat(X_i)\}^2.
\end{equation}
\end{enumerate}
\end{lem}

The weighted decompositions  (\ref{sample-variance-decomp2}) and (\ref{sample-pred-error-decomp2}) in the above lemma hold for any set of nonnegative weights $w_1,\ldots, w_n$ satisfying $\sum_{i=1}^n w_i =1$.  The next lemma shows that  for a particular set of weights defined by (\ref{weight}) below, the decompositions (\ref{sample-variance-decomp2})
and (\ref{sample-pred-error-decomp2}) can be viewed as right-censored data analogs of  the variance decomposition (\ref{sample-variance-decomp}) and  the prediction error decomposition (\ref{sample-pred-error-decomp}), respectively.

\begin{lem}
\label{lem4}
Let
\begin{equation}
\label{weight}
w_i = \frac{\frac{\delta_i}{\hat{G}(T_i-)}}{\sum_{j=1}^n  \frac{\delta_j}{\hat{G}(T_j-)}}, \quad i=1,...,n,
\end{equation}
where $\hat{G}$ is the Kaplan-Meier \citep{kaplan1958} estimate of $G(c)=P(C>c)$.
Assume (\ref{C1}) and (\ref{C2})  hold. Assume further that $C$ is independent of $X$.
Then, under mild regularity conditions,
\begin{eqnarray*}
& &\sum_{i=1}^n w_i\{T_i -\yw\}^2  \convinprob  var(Y);\\
& &\sum_{i=1}^n w_i \{\cmhatw(X_i)- \yw\}^2 \convinprob E\{\cp(X) -\mu_Y\}^2;\\
& &\sum_{i=1}^n w_i \{T_i - \cmhatw (X_i)\}^2 \convinprob E\{Y-\cp(X) \}^2; \\
& & \sum_{i=1}^n w_i\{T_i - \mhat (X_i) \}^2  \convinprob E\{Y-\P(X) \}^2;\\
& & \sum_{i=1}^n w_i \{\cmhatw(X_i)- \mhat(X_i)\}^2  \convinprob E\{\cp(X)-\P(X) \}^2.
\end{eqnarray*}
\end{lem}

Motivated by Lemmas \ref{prop5} and \ref{lem4}, we define the following prediction accuracy measures of $\P(X)$ for right-censored data.
\begin{definition}
\label{summary-censored}
The right censored sample version of $\rho^2$ and $\lambda^2$ are defined by
\begin{equation}
\label{Rsquare2}
R^2_{\mhat} =  \frac{\sum_{i=1}^n w_i \{\cmhatw(X_i)- \yw\}^2}{\sum_{i=1}^nw_i\{T_i -\yw\}^2},
\end{equation}
and
\begin{equation}
\label{Lsquare2}
L^2_{\mhat}=\frac{\sum_{i=1}^n w_i\{T_i -\cmhatw(X_i)\}^2}{\sum_{i=1}^n w_i\{T_i - \mhat (X_i) \}^2},
\end{equation}
where the weight $w_i$'s are defined by (\ref{weight}) and $\cmhatw$ is defined by  (\ref{Q1}). The above defined measures are
interpreted as  the proportion of sample variance of $Y$ explained by $\cmhatw(X)$ and  the proportion of sample mean squared prediction error of $\mhat(X)$ explained by $\cmhatw(X)$, respectively.
\end{definition}

By definition,
$0\le R^2_{\mhat} \le 1$ and $0\le L^2_{\mhat} \le 1$.

\begin{theo}
\label{asymptotics-censored}

\begin{enumerate}
\item[(a)] (Uncensored Data). If there is no censoring, then formulas (\ref{Rsquare2}) and (\ref{Lsquare2})  reduce to  the uncensored data definitions (\ref{Rsquare}) and (\ref{Lsquare}), respectively.

\item[(b)] (Consistency).
Assume the assumptions of Lemma~\ref{lem4}  hold. Then, under mild regularity conditions,
as $n\to\infty$,
$$
R^2_{\mhat}  \convinprob \rho^2_{\P}, \quad\mbox{and}\quad L^2_{\mhat}  \convinprob \lambda^2_{\P}.
$$

\item[(c)]
(Asymptotic normality). Assume the assumptions of Lemma~\ref{lem4} hold.
Then, under certain  regularity conditions,
\begin{eqnarray*}
\sqrt{n} (R^2_{\mhat}  -\rsqp) \convindist N(0, v_{\rho}^2),
\quad\mbox{and}\quad
\sqrt{n} (L^2_{\mhat}  -\lsqp) \convindist N(0, v_{\lambda}^2),
\end{eqnarray*}
as $n\to\infty$, where $v_{\rho}^2$ and $v_{\lambda}^2$ are the asymptotic variances.
\end{enumerate}
\end{theo}

\begin{remark}
\label{remark3.2}
It follows from Theorem \ref{asymptotics-censored} (b) and (c) that the  $R^2_{\mhat}$ and $L^2_{\mhat}$ measures defined by (\ref{Rsquare2}) and (\ref{Lsquare2}) for right censored data are consistent estimates of the population $\rsqp$ and $\lsqp$, respectively, provided that $C$ is independent of $X$ and $Y$. In the next section, we demonstrate by simulation that the
$R^2_{\mhat}$ and $L^2_{\mhat}$ measures are quite robust even if  $C$ depends the covariates.  Furthermore, one could  replace the Kaplan-Meier estimate $\hat{G}(c)$ in (\ref{weight}) by a model-based consistent estimate
$\hat{G}(c|x)$ of $G(c|x)=P(C > c|X=x)$  when there is plausible evidence that $C$ depends on some covariates. In such a case, Theorem  \ref{asymptotics-censored} (b) and (c) would still hold if
$\sup_{c,x} | \hat{G}(c|x) - G(c|x)| \convinprob 0 \quad \mbox{as $n\to\infty$}.$
\end{remark}

\comment{
\begin{remark}
\label{remark3.3}
{\bf We note that condition (\ref{param-assumption1}) holds for most commonly used survival regression models including parametric AFT models, Cox's regression model, Aalen's additive risks model, semiparametric additive risks models (lin and Ying, McKeague), and transformation models (reference).}

\end{remark}
}

\section{Simulations}
\label{simulation}

 In the first simulation, we examine the prediction power of a Cox model by simulating its population $\rho^2_{NP}$ value as defined by (\ref{nonparametricR}) and use it as a benchmark to evaluate the performance of two existing $R^2$-type measures  proposed by
\citet{schemper2000predictive} and \citet{stare2011measure} under
a variety of Cox's models. Specifically, the event time $Y$  is generated from a Cox proportional hazard model:$Y=H_0^{-1}[-\log(U)\times \exp(-\beta^T X)] $, where  $U \sim U(0,1)$, $H^{-1}_0(t)=2t^{1 \over \nu}$ is the inverse function of a Weibull cumulative  hazard function $H_0(t) = (0.5 t)^{\nu}$, and  $X$ is dichotomous = 10* Bernoulli(0.5). We consider six settings by varying $\beta=0.2,  5$, and $\nu=0.5$ (models 1 and 4), 1 (models 2 and 5),  and 10 (models 3 and 6). We  approximate  an population $\rho^2$ value  by averaging  its sample $R^2$ values over 100 Monte Carlo samples of size
$n=5,000$  with no censoring. The results are summarized in Table 1.

\comment{}

\begin{center}
[Insert Table 1  approximately here]
\end{center}

\comment{
\begin{table}[ht]
\label{table1}
\begin{center}
\caption{Simulated Population Proportion $\rho^2_{NP}$ of Explained Variance by the Cox (1972) Model   and the Population Values of $R^2_{SPH}$ and $R^2_{SH}$ Proposed by \citet{schemper2000predictive} and \citet{stare2011measure}.}
 \begin{tabular}{ccccc}
    \hline\hline
  Model &  $\beta$  & $\rho_{NP}^2$ & $R_{SPH}^2$ & $R^2_{SH}$ \\
    \hline
  1&  0.2  & 0.089 & 0.380  & 0.275  \\
   2& 0.2  &  0.271 &  0.381 &  0.276 \\
   3 &0.2  &  0.407 & 0.381  &  0.276 \\
    \cline{2-5}
 4&   5  & 0.091  &  0.499 & 0.502  \\
  5&  5  &  0.332 & 0.500  & 0.505   \\
  6&  5 &  0.971 & 0.500  & 0.503  \\
    \hline
  \end{tabular}\\
  \end{center}
  \end{table}
}

It is seen from Table 1 that the predictive power of a Cox model  depends not only on the regression coefficient $\beta$ (or hazard ratio $e^{\beta}$),  but also on its baseline hazard $h_0(t)$.  A larger $\beta$  does not
always imply a larger proportion of explained variance  when the models are not nested with different baseline hazards (Model 4 versus Model 3). Table 1 also reveals that the $R^2$-type measures proposed
by \citet{schemper2000predictive} and \citet{stare2011measure}
are not effective measures for comparing unnested Cox models. For example, they both are unable to distinguish  between models 4, 5 and 6 as the true proportion
$\rho_{NP}^2$ of explained variance
ranges from 0.09 to 0.97.

In the second simulation, we consider a model with independent censoring to investigate the performance of our proposed sample
prediction accuracy measures $R^2$ and $L^2$
for right-censored data in comparison with the pseudo $R^2$ measures proposed by  \citet{schemper2000predictive} and \citet{stare2011measure} using the population
$\rho^2$ and $\lambda^2$ as benchmarks. Specifically, the event time $Y$ is generated from a Weibull model $\log(Y)=\beta^T X +\sigma  W$, where $\beta=1$, $\sigma=0.15$, $X \sim U(0,1)$, and $W$ has the standard extreme value distribution.  Independent right-censoring time is set to be $C \sim Weibull (shape = 1, scale= b)$. We adjust $b$ to produce censoring rates $25\%$, $0\%$, $50\%$ and $70\%$. We then compute prediction accuracy measures for the Cox PH model that is well specified and for the log-normal AFT model that is obviously mis-specified. Again, the population
$\rho^2$ and $\lambda^2$  are approximated by the averaged sample values over 100 Monte Carlo samples of size $n=5,000$, assuming no censoring. For the sample measures,  we consider sample size n = (50, 200, 500) for each of the parameter settings. The results are reported in Table 2. Each entry in Table 2 is based on 1,000 replications.


\begin{center}
[Insert Table 2 approximately here]
\end{center}
First, we observe from Table 2  that the sample $L^2$ and $R^2$ measures for both censored and uncensored
data estimate the corresponding population values well with  small bias  across almost all scenarios considered except
when there is heavy censoring. Secondly, $L^2$ effectively captures the facts that the Cox model is correctly specified
($L^2=1$)
and that the log-normal AFT model is mis-specified and the predictor is not the mean response ($L^2= 0.789$). Finally, the $R^2$ measures proposed by \citet{schemper2000predictive} and \citet{stare2011measure}
do not really measure the proportion of explained variance, which is consistent with what is observed from the previous simulation
(Table 1).  In particular, the measure $R^2_{SPH}$ of  \citet{schemper2000predictive}  has the same value  for the Cox model and the log-normal
AFT model  and thus is unable to distinguish between the prediction power of these two models.

In the third simulation, we study the robustness of the $R^2$ and $L^2$ measures defined in Section 3 when the independent
censoring assumption is perturbed. The  simulation  setup is similar to the second simulation except that the censoring time $C$ is dependent on the covariate $X$ and that $Y$ and $C$ are conditional independent given the covariate. Specifically, $\log(C)=\gamma^T_c X +\theta_c \times V$, where $X \sim U(0,1)$, $\theta_c=4$, $V \sim$ extreme value distribution, and $\gamma_c$ is adjusted to give censoring rates $25\%$, $50\%$ and $70\%$.  The results are presented in Table 3.


\begin{center}
[Insert Table 3 approximately here]
\end{center}

It is seen that the results in Table 3  are very similar to Table 2. Therefore our proposed $R^2$ and $L^2$ measures 
are not very sensitive to violations of the independent censoring assumption.

Finally, we also conducted simulations when the Kaplan-Meier estimate $\hat{G}$ in (\ref{weight}) is replaced by a Cox model
based estimate of the conditional survival function of $C$. The results are similar and thus not included here.

\section{Real Data Examples}

\begin{ex}
\label{MooreLaw}
{\rm
(Moore's Law). Moore's law predicts that the number of transistors in a dense integrated circuit doubles approximately every two years \citep{moore1975progress,schaller1997moore}. A scatter plot of the $\log_2$-transformed
transistor count together with the fitted least squares line from year 1971 to 2012 is depicted in Figure 2(a). The
$R^2$ for the linear model prediction of the $\log_2$-transformed
transistor count is 0.98, such that 98\% of the variation in the $\log_2$-transformed
transistor count is explained by the fitted least squares line. The corresponding $L^2$ is 1 as expected for a linear model.
In contrast, if one is interested in the prediction  of the untransformed transistor count, then $R^2=0.69$ (Figure 2(b)), meaning that only 69\% of the variation in the untransformed transistor count is explained by the power prediction function {\tt $Y=2^{a +b x}$ } after a linear correction. The log-linear model for the untransformed transistor count has an $L^2=0.96$, so that the linear correction makes very little improvement over the uncorrected prediction.

\begin{center}
[Insert Figure 2 approximately here]
\comment{
\begin{figure}[htbp]
\begin{center}
\includegraphics[scale=0.60]{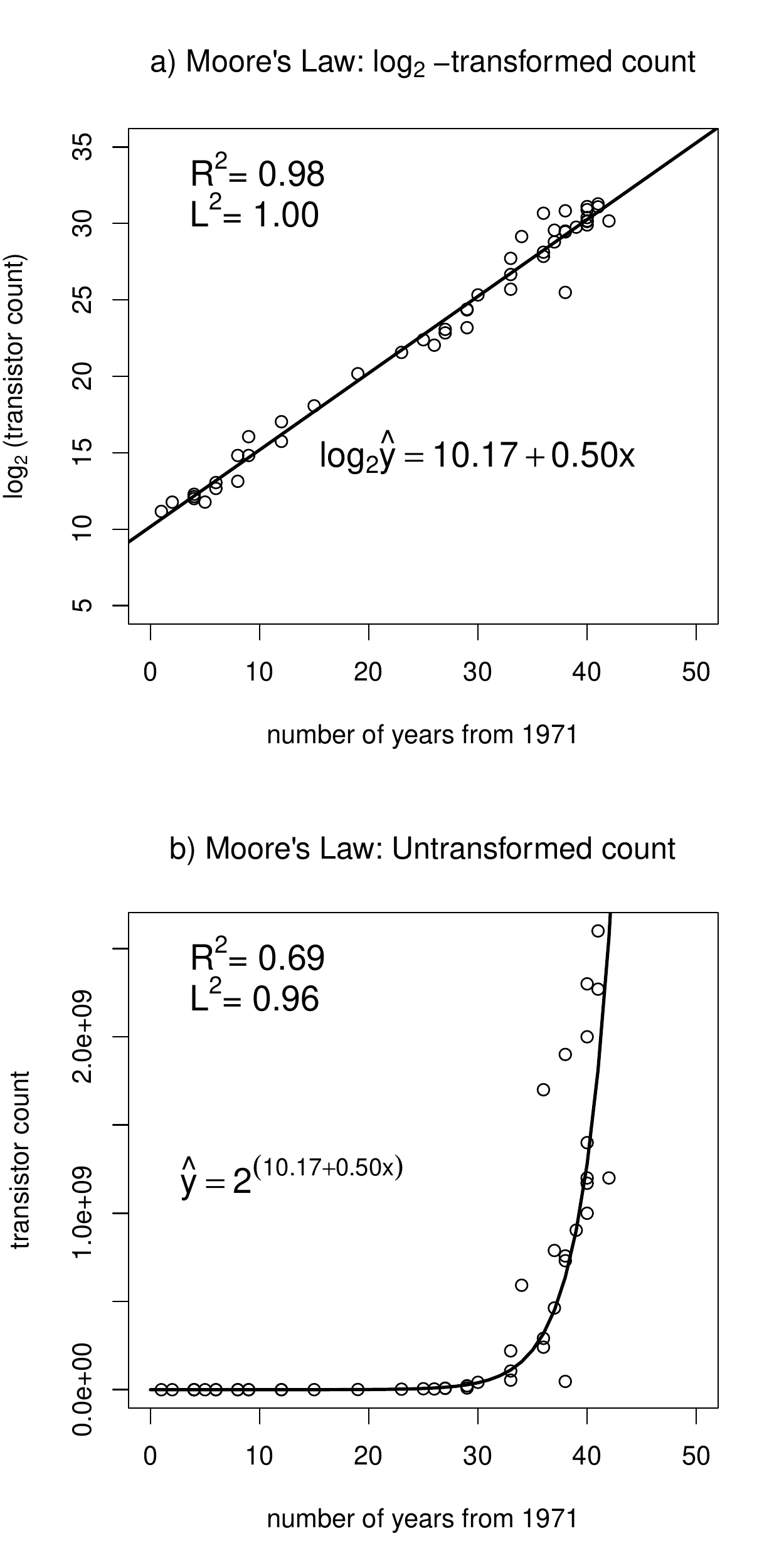}
\caption{\label{fig-example1a} (Moore's Law data)  (a) Prediction power of the log-transformed $Y$; (b) Prediction power of the untransformed $Y$}
\end{center}
\end{figure}
}
\end{center}

}
\end{ex}

\begin{ex}
\label{OvarianCancer}
{\rm
(NY-ESO-1 for Ovarian Cancer) The cancer testis antigen NY-ESO-1 is a potential target for cancer immunotherapy
and has been the focus of multiple cancer vaccine studies. An important question is whether NY-ESO-1 is an important
prognostic marker for overall survival. Table 4 presents the Cox regression results of overall survival
based on a right-censored data from 36 platinum resistant ovarian cancer patients treated at UCLA.
\begin{center}
[Insert Table 4 approximately here]
\end{center}
It is seen from Table 4 that NY-ESO-1 is statistically significant (p-value=0.04) at an $\alpha=0.05$ level with a hazard ratio 3.12. However, as demonstrated in Section 4 (Table 1), a large hazard ratio does not always imply high prediction power. To
evaluate the prediction power of NY-ESO-1 on overall survival, we computed the prediction accuracy measures $R^2$ and $L^2$
of two Cox's models  with and without NY-ESO-1 in Table 5, which shows that the $R^2$ value drops from 0.48 to 0.36 when NY-ESO-1 is removed from the model, indicating NY-ESO-1 is a potentially important prognostic marker for overall survival.

\begin{center}
[Insert Table 5 approximately here]
\end{center}

We also investigated if CA 125, a protein tumor marker measured in the blood, is a good prognostic marker for
overall survival of the same patient population. By comparing models, with and without CA 125, we see that
the $R^2$ value drops only minimally from  0.483 to 0.477 when CA 125 is removed from the model. Hence, there is  no evidence of CA 125
being a good  prognostic marker for overall survival even though it has a larger hazard
ratio (3.92) than that (3.12) of NY-ESO-1, which is not surprising for unnested Cox's models with different baseline hazards as observed in Section 4  (Table 1) . We also note that the $L^2$ values for the Cox models are all 96\%, or higher,  indicating that there is little or no need for a linear correction.

\comment{
It is worth noting that we would not reach the same conclusions when the pseudo $R^2$ values $R^2_{SH}$ of \citet{schemper2000predictive} and $R^2_{SPH}$ of \citet{stare2011measure} were used, which is not surprising given that they demonstrated some undesirable properties in our simulations.  In fact, they gave almost the opposite conclusions. For example,
$R^2_{SH}$ drops only minimally from 0.328 to 0.295 if NY-ESO-1 is removed from the model and but substantially from 0.328 to 0.218 if CA 125 is instead removed from the model.  Similarly, $R^2_{SPH}$ drops only minimally from 0.515 to 0.473 if NY-ESO-1 is removed from the model and
substantially from 0.515 to 0.396 if CA 125 is instead removed from the model. {\tt Do we need to include this paragraph at all since we already have the simulation to indicate that these measures are not good?}

In addition, Table 4 (panel B) also presents similar analysis using the log-normal accelerated failure time model
$log(Y) = \alpha +\beta^T X +\epsilon$ where $\epsilon\sim N(0,\sigma^2)$ and the prediction function
$\exp\{ \hat\alpha + \hat\beta X\}$. {\tt not sure if we want to include this part since it would be difficult to
compare two models when their $L^2$ values are different}
}
}
\end{ex}

\begin{ex}
\label{PBC}
(Comparison of Feature Selection Methods).
{\rm
In this example, we use the right censored primary biliary cirrhosis (PBC) data \citep{tibshirani1997lasso, therneau2000modeling} to illustrate how the proposed prediction accuracy measures can be used  to compare different feature selection methods for high dimensional data. The PBC data is from the Mayo Clinic trial in primary biliary cirrhosis  of the liver conducted between 1974 and 1984. Similar to  \citet{tibshirani1997lasso}, we use 276 patients after removing missing observations. We consider  153 features that include 17 main effects and 136 two-way interactions. Table 6 summarizes the prediction accuracy statistics of models selected by three popular feature selection methods for the Cox model: LASSO \citep{tibshirani1997lasso},
SCAD \citep{fan2002variable}, and Adaptive LASSO \citep{zhang2007adaptive}. 

\begin{center}
[Insert Table 6 approximately here]
\end{center}

It is seen from Table 6 that with a linear correction, the model selected by Adaptive LASSO  uses the fewest (13)  features to achieve the highest proportion of explained variation  ($R_{A-LASSO}^2=0.50$). In contrast, the model selected by
 LASSO  uses 11 more features to achieve a slightly lower  $R_{LASSO}^2=0.49$.  The linear correction is needed for the Adaptive LASSO model ($L^2_{A-LASSO}=0.84$), but does not  seem to be necessary for the LASSO model ($L_{LASSO}^2=0.94$). The model selected by SCAD  is the least desirable in this example since it has the lowest $R_{SCAD}^2=0.45$ and $L^2_{SCAD}=0.77$.

\comment{
\begin{table}[ht]
\label{table6}
\centering
 \caption{Prediction accuracy measures for three Cox's models selected using LASSO,
SCAD, and Adaptive LASSO, respectively, for the primary biliary cirrhosis (PBC) data}
\scalebox{0.85}{
\begin{tabular}{cccc}
  \hline\hline
  & $\# $ of Selected Features  & $R^2$& $L^2$  \\
\hline
  LASSO & 24 & 0.49&  0.94  \\
       \hline
 SCAD & 14 & 0.45 &0.77\\ \hline
 Adaptive LASSO&13 &0.50 &0.84\\
       \hline
      \end{tabular}
}
\end{table}
}
}
\end{ex}

\section{Discussion}
We have introduced a pair of accuracy measures for the predictive power of a prediction function
based on a possibly mis-specified regression model.   Both population and  sample measures are derived. The first measure $\rho^2$ is an extension of the classical $R^2$ statistic for a linear model, quantifying the amount of variability in the response that is explained by a linearly corrected prediction  function. The second measure $\lambda^2$ is the proportion of the squared prediction error of the original prediction function that is explained by the corrected prediction function, quantifying the distance between the corrected and uncorrected predictions. Generally speaking,  $\rho^2$ measures the prediction function's ability to capture the variability of the response and $\lambda^2$ measure its bias  for predicting the mean regression function.  When used together, they give a
complete accuracy of the predictive power of a prediction function.

We have also extended the proposed prediction accuracy measures to  right-censored data by deriving right-censored sample versions of  the variance and  prediction error decompositions.  As discussed earlier, the resulting prediction accuracy measures for right-censored data possess many  appealing properties that other existing pseudo $R^2$ measures do not have: 1) for the linear model,  our $R^2$ statistic  reduces to the classical coefficient of determination when there is no censoring; 2) If the prediction is the conditional mean response based on a correctly specified model ,  then our $R^2$ statistic is a consistent estimate of the population nonparametric coefficient of determination or the proportion of variance of $Y$ explained by $E(Y|X)$; 3) our method is applicable to any event time model;  4) our measures are defined without requiring the model to be correctly  specified, and 5)  our measures can be used to compare unnested models. 

We have implemented our methods for right-censored data using {\it R}. Our  {\it R} code is available upon request.

Lastly,  this paper focuses on  i.i.d.\  complete data and  right censored data. Future efforts to develop prediction accuracy measures for correlated data such as longitudinal data  and for other censoring patterns are warranted.

\bigskip
\begin{center}
{\large\bf SUPPLEMENTARY MATERIAL}
\end{center}

\begin{description}

\item[Appendix:] Proofs of the lemmas and theorems. (pdf)
\end{description}
\bibliographystyle{natbib}
\bibliography{reference}


\newpage
\begin{table}[ht]
\label{table1}
\begin{center}
\caption{Simulated Population Proportion $\rho^2_{NP}$ of Explained Variance by the Cox (1972) Model   and the Population Values of $R^2_{SPH}$ and $R^2_{SH}$ Proposed by \citet{schemper2000predictive} and \citet{stare2011measure}.}
 \begin{tabular}{ccccc}
    \hline\hline
  Model &  $\beta$  & $\rho_{NP}^2$ & $R_{SPH}^2$ & $R^2_{SH}$ \\
    \hline
  1&  0.2  & 0.089 & 0.380  & 0.275  \\
   2& 0.2  &  0.271 &  0.381 &  0.276 \\
   3 &0.2  &  0.407 & 0.381  &  0.276 \\
    \cline{2-5}
 4&   5  & 0.091  &  0.499 & 0.502  \\
  5&  5  &  0.332 & 0.500  & 0.505   \\
  6&  5 &  0.971 & 0.500  & 0.503  \\
    \hline
  \end{tabular}\\
  \end{center}
  \end{table}

\newpage
\begin{table}[ht]
\label{table2}
\centering
 \caption{(Independent Censoring) Simulated Prediction Accuracy Measures for the Cox Model and for the Log-Normal Accelerated Failure Time
 (AFT) Model.}
\scalebox{0.85}{
\begin{tabular}{cccccccccc}
  \hline\hline
     & &  \multicolumn{4}{c}{Cox's Model (Correctly Specified)} & &  \multicolumn{3}{c}{Log-normal AFT Model (Mis-specified)} \\
    \cline{3-6}  \cline{8-10}
 CR & N & $L^2$ & $R^2$ & $R^2_{SPH}$ & $R^2_{SH}$ & & $L^2$ & $R^2$ & $R^2_{SPH}$ \\
  \hline
 0\% & $\infty$ & 100.0 & 70.4 & 65.4 & 50.3 && 78.9 & 70.4 & 65.4 \\
 \cline{1-10}
   0\% & 50 & 96.6(1.5) & 70.7(7.4) & 65.2(5.1) & 49.2(5.9) && 75.9(18.6) & 70.6(7.6) & 65.2(5.1) \\
    & 200 & 99.6(0.3) & 70.6(3.9) & 65.4(2.5) & 50.1(3.0) &&  77.8(10.6) & 70.5(3.9) & 65.4(2.5) \\
    & 500 & 99.9(0.1) & 70.5(2.3) & 65.4(1.5) & 50.3(1.8) && 78.2(7.2) & 70.5(2.3) & 65.4(1.5) \\
   25\% & 50 & 96.4(3.0) & 70.6(8.9) & 65.4(6.0) & 47.7(7.3) & & 73.7(20.9) & 70.3(9.0) & 65.4(6.0) \\
    & 200 & 99.5(0.5) & 70.7(4.5) & 65.4(2.7) & 49.8(3.4) &&  76.9(11.9) & 70.6(4.5) & 65.4(2.7) \\
    & 500 & 99.9(0.2) & 70.6(2.7) & 65.4(1.7) & 50.2(2.1) & & 77.7(8.3) & 70.6(2.7) & 65.4(1.7) \\
  50\% & 50 & 93.5(5.9) & 71.4(11.0) & 66.0(7.6) & 47.8(8.6)  &&  69.2(24.9) & 70.9(11.2) & 66.0(7.6) \\
    & 200 & 99.0(1.1) & 70.8( 5.3) & 65.6(3.2) & 49.9(3.8) &&  74.9(15.0) & 70.7(5.4) & 65.6(3.2) \\
  & 500 & 99.7(0.3) & 70.6(3.3) & 65.5(2.0) & 50.1(2.4)  && 76.5(9.9) & 70.6(3.3) & 65.5(2.0) \\
   70\% & 50 & 87.7(12.7) & 69.2(15.3) & 65.9(10.1) & 45.9(11.1) & &58.6(27.8) & 68.3(15.8) & 65.9(10.1) \\
     & 200 & 97.5(3.5) & 70.5(7.2) & 65.6(4.3) & 49.2(4.8) && 72.3(18.4) & 70.3(7.4) & 65.6(4.3) \\
    & 500 & 99.3(0.9) & 70.8(4.5) & 65.6(2.6) & 50.2(3.0) &&  74.3(13.3) & 70.7(4.5) & 65.6(2.6) \\
\hline
\end{tabular}
}
\end{table}

\newpage
\begin{table}[ht]
\label{table3}
\centering
 \caption{(Dependent Censoring) Simulated Prediction Accuracy Measures for the Cox Model and for the Log-Normal Accelerated Failure Time
 (AFT) Model.}
\scalebox{0.85}{
\begin{tabular}{cccccccccc}
  \hline \hline
    & &  \multicolumn{4}{c}{Cox's Model (Correctly Specified)} & &  \multicolumn{3}{c}{Log-normal AFT Model (Mis-specified)} \\
    \cline{3-6}  \cline{8-10}
 CR & N & $L^2$ & $R^2$ & $R^2_{SPH}$ & $R^2_{SH}$ & & $L^2$ & $R^2$ & $R^2_{SPH}$ \\
 \hline
0\% & $\infty$ & 100.0 & 70.4 & 65.4 & 50.3 &   & 78.9 & 70.4 & 65.4 \\
\cline{1-10}
  25\% & 50 & 96.5(2.2) & 68.0(9.1) & 63.6(6.4) & 49.8(6.7) &   & 71.6(23.4) & 67.7(9.4) & 63.5(7.6) \\
   & 200 & 99.6(0.4) & 67.5(4.6) & 63.6(3.0) & 50.6(3.4) &   & 75.0(16.6) & 67.3(5.1) & 63.5(5.0) \\
   & 500 & 99.9(0.1) & 67.6(2.9) & 63.7(1.8) & 50.8(2.1) &   & 76.9(11.0) & 67.6(2.9) & 63.7(1.8) \\
   50\% & 50 & 93.5(4.7) & 69.9(11.2) & 64.7(8.0) & 50.2(8.4) &   & 70.3(25.6) & 69.3(11.7) & 64.3(10.6) \\
  & 200 & 99.3(0.8) & 69.3(5.4) & 64.8(3.5) & 51.1(3.9) &   & 76.2(16.2) & 69.0(6.6) & 64.4(7.8) \\
   & 500 & 99.8(0.2) & 68.9(3.3) & 64.6(2.1) & 51.0(2.4) &   & 76.9(11.5) & 68.6(5.9) & 63.8(10.3) \\
  70\% & 50 & 84.0(12.8) & 71.0(15.4) & 65.1(11.7) & 48.4(12.5) &   & 65.4(27.2) & 70.2(15.4) & 65.1(11.7) \\
    & 200 & 98.2(1.7) & 71.0(7.0) & 65.4(4.5) & 49.9(5.3) &   & 75.0(17.7) & 70.8(7.1) & 65.4(4.5) \\
    & 500 & 99.5(0.5) & 70.8(4.3) & 65.4(2.7) & 50.1(3.2) &   & 76.7(12.2) & 70.7(4.3) & 65.4(2.7) \\
\hline
\end{tabular}
}
\end{table}

\newpage
\begin{table}[ht]
\label{table4}
\centering
 \caption{Cox's proportional hazards regression of overall survival
based on a right-censored data from 36 platinum resistant ovarian cancer patients treated at UCLA}
\scalebox{0.85}{
\begin{tabular}{ccccccccc}
  \hline\hline
 & \multicolumn{2}{c}{Full Model} & &\multicolumn{2}{c}{Reduced Model} && \multicolumn{2}{c}{Reduced Model}\\
  & \multicolumn{2}{c}{} & &\multicolumn{2}{c}{Without NY-ESO-1} && \multicolumn{2}{c}{Without CA 125}\\
\cline{2-3} \cline{5-6}\cline{8-9}
variables                & HR & p-value && HR & p-value && HR & p-value\\
\hline
  stage($3\&4$ vs $1\&2$) & 4.45 & 0.10 && 7.86& 0.02  & &3.97& 0.10\\
  grade($1\&2$ vs 3)      & 1.07 &  0.89 && 1.00 & 0.99 & &0.86& 0.76\\
  histology &   &  &   &  & & \\
   $\, \, \,$ endometrioid vs clear cell & 0.95 & 0.95&& 0.42 & 0.28 && 1.34& 0.72\\
   $\, \, \,$  serious vs clear cell & 0.29 & 0.09 && 0.21 & 0.04&& 0.58& 0.41 \\
   preop CA125 ($>500$ vs $\le500$) & 3.92 & 0.01 && 4.17 & $<$0.01 &&-- & --\\
   NY-ESO1 ($>12$ vs $\le12$) & 3.12 & 0.04 && -- & -- & &3.67 & 0.02 \\
  \hline
\end{tabular}
}
\end{table}

\newpage
\begin{table}[ht]
\label{table5}
\centering
 \caption{Prediction accuracy measures for Cox's proportional hazards models 
based on a right-censored   ovarian cancer data }
\scalebox{0.85}{
\begin{tabular}{cccccc}
  \hline\hline
  & $R^2$  & $L^2$ \\
\hline
Full Cox's Model WIth All Variables      &  0.483 & 0.991	 \\
       \hline
 Reduced Cox's Model  Without NY-ESO-1    &  0.363	&0.991\\
         \hline
  Reduced Cox's Model Without CA 125      & 0.477 & 0.963 \\
       \hline
      \end{tabular}
}
\end{table}

\newpage
\begin{table}[ht]
\label{table6}
\centering
 \caption{Prediction accuracy measures for three Cox's models selected using LASSO,
SCAD, and Adaptive LASSO, respectively, for the primary biliary cirrhosis (PBC) data}
\scalebox{0.85}{
\begin{tabular}{cccc}
  \hline\hline
  & $\# $ of Selected Features  & $R^2$& $L^2$  \\
\hline
  LASSO & 24 & 0.49&  0.94  \\
       \hline
 SCAD & 14 & 0.45 &0.77\\ \hline
 Adaptive LASSO&13 &0.50 &0.84\\
       \hline
      \end{tabular}
}
\end{table}

\newpage
\begin{center}
\begin{figure}[htbp]
\begin{center}
\includegraphics[scale=0.65]{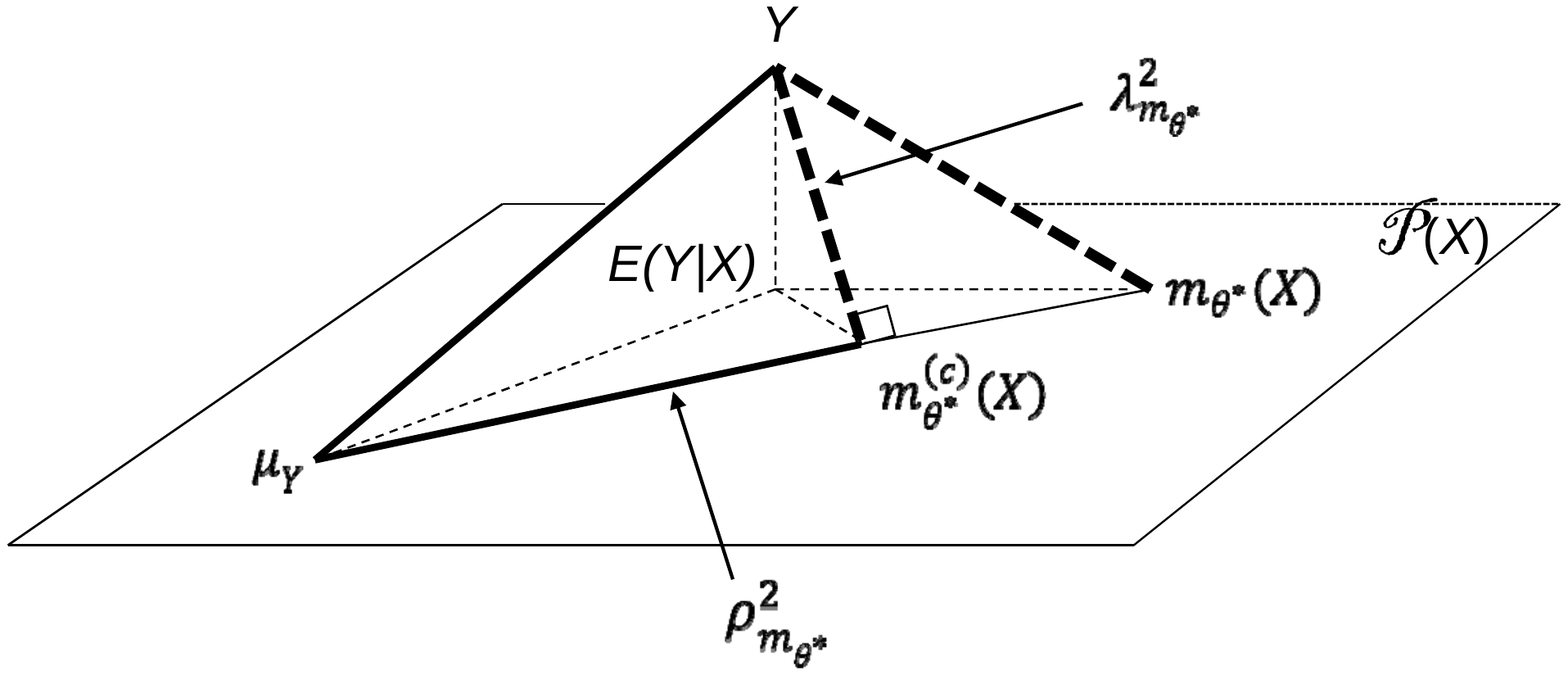}
\caption{\label{fig1} Geometric interpretation of $\rsqp$ and $\lambda^2_{\P}$}
\end{center}
\end{figure}
\end{center}

\newpage
\begin{center}
\begin{figure}[htbp]
\begin{center}
\includegraphics[scale=0.60]{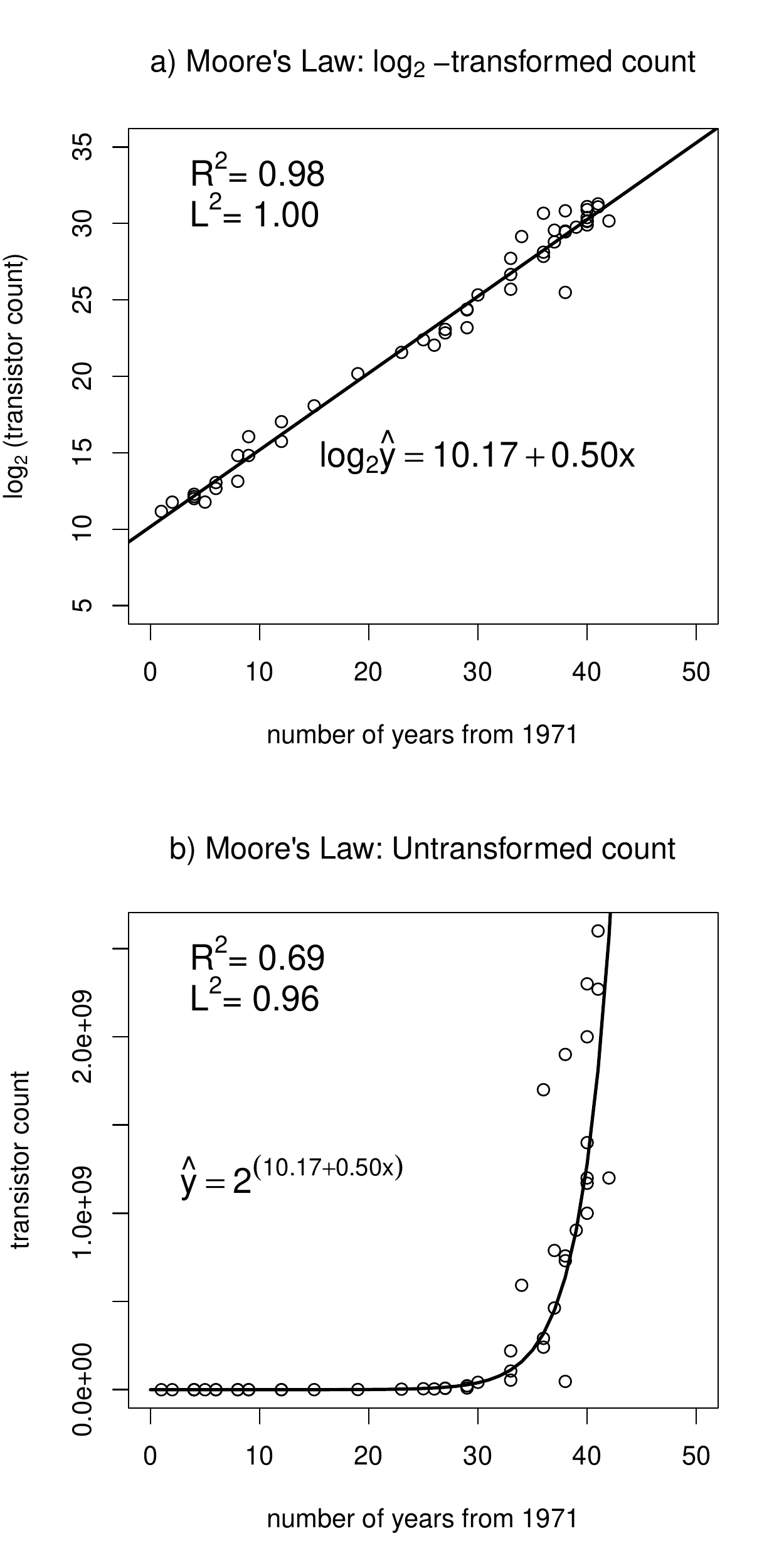}
\caption{\label{fig-example1a} (Moore's Law data)  (a) Prediction power of the log-transformed $Y$; (b) Prediction power of the untransformed $Y$}
\end{center}
\end{figure}
\end{center}

\newpage

\cleardoublepage
\appendix

\numberwithin{equation}{section}
\numberwithin{figure}{section}

\makeatletter   
 \renewcommand{\@seccntformat}[1]{APPENDIX~{\csname the#1\endcsname}.\hspace*{1em}}
\makeatother


\section{Supplementary Material}
\pagenumbering{arabic}

PROOF OF LEMMA \ref{prop1}. (a) Note that
\begin{eqnarray*}
var(Y)& =&E(Y-\mu_Y)^2\\
&=& E\{Y-\cp(X)\}^2  + 2 E\{ \cp(X) -\mu_Y\} \{ Y -\cp(X)\} + E\{ \cp(X) -\mu_Y\}^2.
\end{eqnarray*}
So it suffices to show that
\begin{equation}
\label{A1}
E\{ \cp(X) -\mu_Y\} \{ Y -\cp(X)\} =0.
\end{equation}
Recall that $\cp(X) = \ta + \tb\P(X)$, where
$(\tilde{a}, \tilde{b})=\arg \min_{\alpha, \beta} E\{ Y -  (\alpha +\beta \P(\X))\}^2 $.  Thus,
\begin{eqnarray*}
\frac{ \partial  E\{ Y -  (\alpha +\beta \P(\X))\}^2 }{ \partial \alpha} \bigg{|}_{(\alpha,\beta)=(\ta, \tb)} =
-2 E\{ Y -  (\ta +\tb \P(\X))\}=0,
\end{eqnarray*}
and
\begin{eqnarray*}
\frac{ \partial  E\{ Y -  (\alpha +\beta \P(\X))\}^2 }{ \partial \beta} \bigg{|}_{(\alpha,\beta)=(\ta, \tb)}
=-2 E[\{ Y -  (\ta +\tb \P(\X))\}\P(X)]=0,
\end{eqnarray*}
which imply that
\begin{eqnarray}
\label{A2}
 E\{ Y -\cp(X)\}=0,
\end{eqnarray}
and
\begin{eqnarray}
\label{A3}
E[\{ Y -\cp(X)\} \P(X)]=0.
\end{eqnarray}
Finally,  (\ref{A1}) follows  from (\ref{A2}) and (\ref{A3}). This proves (\ref{decomposition1}).

(b).  Note that
\begin{eqnarray*}
& &E\{ Y-\cp(X)\} \{\cp(X) -\P(X)\}\\
&=& E\{ Y-\cp(X)\} \{\ta + \tb \P(X) -\P(X)\} \\
&=& \ta E\{ Y-\cp(X)\}  + (\tb -1)E[\{ Y-\cp(X)\} \P(X)] \\
&=& 0,
\end{eqnarray*}
where the last equality follows from (\ref{A2}) and (\ref{A3}). This immediately implies that
 (\ref{decomposition2}) holds. $\qquad \Box$

\bigskip
PROOF OF THEOREM \ref{prop2}.
The proofs for parts (a)-(c) are straightforward. Part (d)  follows directly from the fact that $\mu(X)=E(Y|\X)$ is the best prediction function for $Y$ among all functions of  $\X$
in a sense that $E\{Y-\mu(X)\}^2 \le E\{Y-Q(X)\}^2$
for any $p$-variate  function $Q$, and that the equality holds when $Q(X)=\mu(X)$. $\Box$

\bigskip
PROOF OF LEMMA \ref{prop3}.
(a). The variance decomposition (\ref{sample-variance-decomp}) is a trivial consequence of the fact that $\cmhat(X)$ is the fitted value from the simple linear regression of $Y$ on $\mhat(X)$.

(b) Now we prove  the prediction error decomposition (\ref{sample-pred-error-decomp}). For the simple linear regression of $Y$ on a covariate $Z$, it is well known that
\begin{equation}
\label{A4}
\sum_{i=1}^n e_i Z_i=0\quad\mbox{and}\quad \sum_{i=1}^n e_i \hat{y}_i=0,
\end{equation}
where $\hat{y}_i$ is the fitted value and $e_i=Y_i - \hat{y}_i$ is the residual  at $Z_i$, $i=1,\ldots,n$.
In our context, $Z_i=\mhat(X_i)$ and $\hat{y}_i=\cmhat(X_i)$, and thus (\ref{A4}) implies that
$$
\sum_{i=1}^n \{Y_i-\cmhat(X_i)\} \P(X_i)=0\quad\mbox{and}\quad \sum_{i=1}^n \{Y_i-\cmhat(X_i)\} \cmhat(X_i)=0.
$$
Consequently,
\begin{eqnarray*}
\sum_{i=1}^n \{Y_i - \mhat(X_i)\}^2
&=& \sum_{i=1}^n \{Y_i - \cmhat(X_i)\}^2 + \sum_{i=1}^n \{\cmhat(X_i) - \mhat(X_i)\}^2 \\
& & + 2 \sum_{i=1}^n \{Y_i - \cmhat(X_i)\} \{\cmhat(X_i) - \mhat(X_i)\}^2 \\
&=&  \sum_{i=1}^n \{Y_i - \cmhat(X_i)\}^2 + \sum_{i=1}^n \{\cmhat(X_i) - \mhat(X_i)\}^2.
\end{eqnarray*}
This proves  (\ref{sample-pred-error-decomp}).
$\Box$

\bigskip
PROOF OF THEOREM \ref{asymptotics}.
 (a) It suffices to show that
\begin{eqnarray}
\label{eq:A5}
\frac{1}{n} \sum_{i=1}^n Y_i \mhat (X_i) \convinprob E\{Y \P(X)\},\\
\frac{1}{n} \sum_{i=1}^n  \mhat (X_i) \convinprob E\{ \P(X)\},\\
\frac{1}{n} \sum_{i=1}^n  \mhat ^2(X_i) \convinprob E\{ \P^2(X)\}.
\end{eqnarray}
We only prove (\ref{eq:A5}) here because the proof of the other two results are similar. Note that
\begin{eqnarray*}
\frac{1}{n} \sum_{i=1}^n Y_i \mhat (X_i)  &=& \frac{1}{n} \sum_{i=1}^n Y_i \P(X_i)
+ \frac{1}{n} \sum_{i=1}^n Y_i \{\mhat (X_i)  -\P(X_i)\} \\
&=& I_1 + I_2.
\end{eqnarray*}
By the law of large numbers, $I_1\convinprob E\{Y \P(X)\}$. Moreover,  under the assumption (\ref{C2}),
\begin{eqnarray*}
|I_2| \le \sup_{x}|\mhat(x) - \P(x) | \left(\frac{1}{n} \sum_{i=1}^n |Y_i |\right)
 \convinprob 0,
\end{eqnarray*}
This proves (\ref{eq:A5}).

(b).  Note that
\begin{eqnarray*}
\frac{1}{\sqrt{n}}  \sum_{i=1}^n [Y_i \mhat (X_i)  - E\{Y \P(X) \} ]&=& \frac{1}{\sqrt{n}}
\sum_{i=1}^n [Y_i \P(X_i)  - E\{Y \P(X) \} ]\\
& &+ \frac{1}{\sqrt{n}} \sum_{i=1}^n Y_i \{\mhat (X_i)  -\P(X_i)\} \\
&=& J_1 + J_2.
\end{eqnarray*}
The asymptotic normality of $J_1$ follows from the Central Limit Theorem. Moreover,  under the assumption (\ref{C2}),
\begin{eqnarray*}
|J_2| \le \sup_{x}|\mhat(x) - \P(x) | \left(\frac{1}{\sqrt{n}} \sum_{i=1}^n |Y_i |\right)
 \convinprob 0.
\end{eqnarray*}
One can indeed establish the joint convergence to a multivariate normal limit of multiple
quantities in the expression of $R^2_{\mhat}$ and $L^2_{\mhat}$. Then part (b) follows from the delta method.
$\quad \Box$

\bigskip
PROOF OF LEMMA \ref{prop5}.
(a) Recall that $W=diag(w_1, \ldots, w_n)$.
Define $\by=(T_1,\ldots,T_n)^{'}$, $\hat{\by}=(\cmhatw(X_1), \ldots, \cmhatw(X_n))^{'}$,
$\bz=(\mhat(X_1), \ldots, \mhat(X_n))^{'}$, and $\bZ=(\b1, \bz)$.
where $\b1=(1,\ldots,1)^{'}$ is a $n$ dimensional column vector of 1's. Then, by the definition of
$\cmhatw$, we have
$$
\hat{\by} = \bZ (\bZ^{'}W \bZ)^{-1} \bZ^{'} W\by.
$$
Note that
\begin{eqnarray*}
(\by\! -\! \hat{\by})^{'} W (\b1\  \bz)=(\by\! -\! \hat{\by})^{'} W \bZ ={\by}^{'} \{I-W\bZ  (\bZ^{'}W \bZ)^{-1} {\bZ}^{'} \}W\bZ= 0,
\end{eqnarray*}
which implies that
\begin{eqnarray}
\label{A5}
(\by - \hat{\by})^{'} W \b1 =0,  (\by - \hat{\by})^{'} W \bz =0,\;
\mbox{and}
(\by\! -\!  \hat{\by})^{'} W \hat{\by} =(\by\! -\! \hat{\by})^{'} W \bZ (\bZ^{'}W \bZ)^{-1} \bZ^{'} W\by= 0.
\end{eqnarray}
Therefore,
\begin{eqnarray*}
 \sum_{i=1}^n w_i\{T_i -\yw\}^2
 &=&(\by - \b1 \b1^{'}W\by)^{'} W (\by - \b1 \b1^{'}W\by)\\
 &=&(\by -\hat{\by})^{'} W (\by -\hat{\by}) + (\hat{\by} -  \b1 \b1^{'}W\by)^{'} W (\hat{\by} -  \b1 \b1^{'}W\by) \\
  & & + 2(\by -\hat{\by})^{'} W (\hat{\by} -  \b1 \b1^{'}W\by)\\
 &=& (\by -\hat{\by})^{'} W (\by -\hat{\by}) + (\hat{\by} -  \b1 \b1^{'}W\by)^{'} W (\hat{\by} -  \b1 \b1^{'}W\by) \\
 &=& \sum_{i=1}^n w_i \{\cmhatw(X_i)- \yw\}^2
+ \sum_{i=1}^n w_i \{T_i - \cmhatw (X_i)\}^2,
\end{eqnarray*}
where the third equality follows from (\ref{A5}). This proves part (a).

(b).
\begin{eqnarray*}
 \sum_{i=1}^n w_i\{T_i-\mhat(X_i)\}^2
 &=&  \sum_{i=1}^n w_i\{T_i -\cmhatw(X_i)\}^2
+\sum_{i=1}^n w_i \{\cmhatw(X_i)- \mhat(X_i)\}^2 \\
    & & + 2\sum_{i=1}^n w_i \{T_i -\cmhatw(X_i)\} \{\cmhatw(X_i)- \mhat(X_i)\}\\
 &=& \sum_{i=1}^n w_i\{T_i -\cmhatw(X_i)\}^2
+\sum_{i=1}^n w_i \{\cmhatw(X_i)- \mhat(X_i)\}^2 \\
&&+ 2  (\by -\hat{\by})^{'} W ( \hat{\by} - \bz) \\
&=&  \sum_{i=1}^n w_i\{T_i -\cmhatw(X_i)\}^2
+\sum_{i=1}^n w_i \{\cmhatw(X_i)- \mhat(X_i)\}^2 ,
\end{eqnarray*}
where the last equality follows from (\ref{A5}). This proves part (b).
$\quad \Box$

\bigskip
PROOF OF LEMMA \ref{lem4}. We first prove the first result  of Lemma \ref{lem4}.
Note that for any function $h(T,X)$ of $(T,X)$,  we have
\begin{eqnarray*}
E\left\{ \frac{\delta h(T,X)}{1-G(T|X)}\right\} &=& E \left[E\left\{ \frac{\delta h(T,X)}{1-G(T|X)}\Big{|}X,Y\right\} \right]\\
&=& E \left[ E\left\{ \frac{\delta h(Y,X)}{1-G(Y|X)}\Big{|}X,Y\right\} \right]\\
&=& E \left\{ \frac{h(Y,X)}{1-G(Y|X)} E( \delta |X,Y) \right\}\\
&=& E \left\{ \frac{h(Y,X)}{1-G(Y|X)} P( C>Y |X,Y) \right\}\\
&=& E \left\{ \frac{h(Y,X)}{1-G(Y|X)} \{1-G(Y|X)\} \right\}\\
&=& E\left\{h(Y,X)\right\}.
\end{eqnarray*}
In particular,   $h(T,X)=1$, $h(T,X)=T$ and $h(T,X)=T^2$, correspond to
\begin{eqnarray*}
E\left\{ \frac{\delta }{1-G(T|X)}\right\} =1, \quad E\left\{ \frac{\delta T }{1-G(T|X)}\right\} =E(Y),
\quad\mbox{and}\quad E\left\{ \frac{\delta T^2 }{1-G(T|X)}\right\} =E(Y^2),
\end{eqnarray*}
which imply that
$
\yw = \sum_{i=1}^n w_i T_i = \frac{\sum_{i=1}^n \frac{\delta_iT_i}{\hat{G}(T_i-0|X_i)}}{\sum_{i=1}^n \frac{\delta_i}{\hat{G}(T_i-0|X_i)}}
\convinprob E(Y),
$
and  $ \sum_{i=1}^n w_i T_i^2 \convinprob E(Y^2) $. Thus,
$$\sum_{i=1}^n w_i\{T_i -\yw\}^2 = \sum_{i=1}^n w_i T_i^2 - \{\yw \}^2 \convinprob E(Y^2) - \{E(Y)\}^2 = var(Y).$$ The proof for the other results of
the lemma are similar and omitted. $\qquad \Box$

\bigskip
PROOF OF THEOREM \ref{asymptotics-censored}.  (a). If there is no censoring, or $\delta_i=1$ for all $i$, then the Kaplan-Meier 
estimate of the survival function of the censoring time
is identical to 1. Thus $w_i=1/n$ for all $i$. The conclusion of (a) follows immediately. 

The proof  of parts (b) and (c) is essentially the same as that of
Theorem \ref{asymptotics}. and thus we omit the details.  $\qquad \Box$

\end{document}